\documentclass[12pt]{article}
\usepackage{amsmath,amssymb,amsthm}

\numberwithin{equation}{section}

\newtheorem{thm}{Theorem}[section]

\newcommand{\n}{\nonumber}

\renewcommand{\a}{\alpha}
\renewcommand{\l}{\lambda}

\renewcommand{\O}{\Omega}
\newcommand{\vb}{\overline{V}}
\newcommand{\Ob}{\overline{\Omega}}
\newcommand{\se}{(SSE)_{\alpha}}
\newcommand{\bb}{\begin{equation}}
\newcommand{\ee}{\end{equation}}
\newcommand{\bq}{\begin{eqnarray}}
\newcommand{\eq}{\end{eqnarray}}
\newcommand{\bqn}{\begin{eqnarray*}}
\newcommand{\eqn}{\end{eqnarray*}}
\begin{document}
\title{ Continuation of the zero set for discretely self-similar solutions to the Euler equations}
\author{Dongho Chae \\
Chung-Ang University\\ Department of Mathematics\\
 Seoul 156-756, Republic of Korea\\
e-mail: dchae@cau.ac.kr}
\date{}
\maketitle

\begin{abstract}
We are concerned on the  study of  the unique continuation type property  for the 3D incompressible Euler equations in the self-similar type form. Discretely self-similar solution is a generalized notion of the self-similar solution, which is equivalent  to a time periodic solution of the time dependent self-similar Euler equations.
We prove  the unique continuation type theorem for the discretely self-similar solutions to the Euler equations in $\Bbb R^3$. More specifically, we suppose there exists an open set $G\subset \Bbb R^3$ containing the origin such that the velocity field $V\in  C_s^1C^{2}_y (\Bbb R^{3+1})$  vanishes on $G\times (0, S_0)$, where $S_0 > 0$ is the temporal period for $V(y,s)$. Then, we show $V(y,s)=0$ for all  $(y,s)\in \Bbb R^{3+1}$. Similar property also holds for the inviscid magnetohydrodynamic system.\\
\ \\
\noindent{\bf AMS Subject Classification Number:} 35Q31, 76B03,
76W05\\
  \noindent{\bf
keywords:}   Euler equations, discretely self-similar solution, continuation of the zero set, inviscid MHD
\end{abstract}
\section{The main theorems}
\setcounter{equation}{0}
\subsection{The Euler equations}
We are concerned on  the Cauchy problem for the
incompressible 3D Euler equations  in $\Bbb R^3$:
$$
(E) \left\{\aligned  &\frac{\partial v}{\partial t }+v\cdot \nabla v =-\nabla p,\label{e1}\\
&\mathrm{div }\, v=0,\\
& v(x,0)=v_0(x),
\endaligned \right.
$$
 where $v(x,t)=(v_1 (x,t), v_2 (x,t), v_3 (x,t))$ is the velocity, $p=p(x,t)$ is the
  pressure, and $v_0(x)$ is the initial data satisfying div$v_0 =0$.  For  (E) it is well-known that  for  the initial data belonging to the standard Sobolev space, $v_0 \in H^m (\Bbb R^3), m>5/2$, the local well-posedness holds(see e.g. \cite{kat}). The question of the finite time singularity for the local classical solution, however,  is still an outstanding open problem(see \cite{bea,con-fef}). For survey books or articles on the study of the finite time  blow-up problem  of (E) we refer e.g. \cite{maj,con, bar}.
For studies of the possibility of  self-similar blow-up or its generalized version, discretely self-similar blow-up for the Euler and the Navier-Stokes equations there are previous studies (e.g. \cite{nec, cha1,cha2, cha3, cha4, cha-shv, cha-tsai, tsai, tsai1}).
Given $(x_*, T)\in \Bbb R^3\times \Bbb R_+$, we consider the  self-similar
transform of (E), which is defined by
given by
\begin{equation}
\label{11}
 v(x, t)=\frac {1}{(T-t)^{\frac{\a}{\a+1}}} V(y,s), \quad
 p(x,t)=\frac {1}{(T-t)^{\frac{2\a}{\a+1}}} P(y,s),
\end{equation}
where
\begin{equation}
\label{12} y = \frac{x}{(T-t)^{\frac{1}{\a+1}}}, \quad s =  \log \left(\frac{T}{T-t}\right),
\end{equation}
and $\alpha \neq -1$.
Substituting (\ref{11})--(\ref{12}) into (E), we obtain the following system  for $(V,P)$:
\begin{equation}\label{} (SSE)_\a \left\{
\aligned
&\frac {\partial V}{\partial s}+ \frac{\a}{\a+1} V +\frac{1}{\a+1}(y \cdot \nabla)V + (V\cdot \nabla )V =-\nabla
P,\\
& \mathrm{div}\,V=0,\\
& V(y,0)=V_0 (y)=T^{\frac{\a}{\a+1}}v_0 (T^{\frac{1}{\a+1}}y).
\endaligned \right.
\end{equation}
The self-similar solution $(v,p)$ of the Euler equations is  a solution of (E) given by (\ref{11})-(\ref{12}), where $(V,P)$ is a {\em stationary} solution of $(SSE)_\a$. In the case $(V,P)$ is nontrivial $T$ is the blow-up time of the solution $(v,p)$, and we say that the solution given by (\ref{11})-(\ref{12}) is a self-similar blowing-up solution.
This is a solution of (E), having  the scale symmetry, namely
\bq\label{13}
\l ^\a v(\lambda x, \lambda ^{\a+1} (T-t))&=&v(x, T-t), \n \\
 \l^{2\a} p(\lambda x, \lambda ^{\a+1} (T-t))&=&p(x, T-t)
\eq
{\em for all} $\l\neq 1, \a \neq -1$.
On the other hand, there exists a generalized notion of the self-similarity for the solution, called discrete self-similarity.
We say a solution $(v,p)$ of (E) is a discretely self-similar solution to (E) if {\em there exists} a pair $(\lambda, \alpha)$ with  $\lambda \neq 1, \a \neq -1$ such that (\ref{13}) holds true.  We find that $(v,p)$ given by  (\ref{11})-(\ref{12}) is a discrete self-similar solution to (E) with $\l \neq 1$ if and only if  $(V,P)$ of $\se$ satisfies the periodicity in time,
$$
V(y, s)=V(y+S_0), \quad P(y,s)=P(y, s+S_0) \quad \forall (y,s)\in \Bbb R^{3+1}
$$
with
\bb\label{period}
 S_0=-(\a+1)\log \l .
\ee
Conversely, for any time periodic solution  $(V,P)$ of $\se$ with the period $S_0\neq 0$ one can generate discretely self-similar solution of (E) with the scaling parameter given by (\ref{period}).
Thus, the question of the existence of nontrivial discretely self-similar solutions is equivalent to that of the existence of nontrivial time-periodic solution to $\se$.
In \cite{cha4} it is shown that the time periodic solution to $\se$ is trivial, namely a constant vector in space variable if we assume
$|\nabla V|\to 0$ as $|y|\to 0$, and $\a <-1$. In the case $\a >-1$  we have the same conclusion of the triviality  of $V$ under the additional decay condition of vorticity, $|\mathrm{curl} \,V|= o(|y|^{-\gamma})$ with $\gamma>\a+1$.
In order to prove this result we used maximum principle in the far field region restricted by a suitable cut-off function.
In this paper we develop the method of \cite{cha4} further to prove the unique continuation type property  for the time periodic solution of $\se$ as follows.
\begin{thm}
Let $\a  \neq -1 $ and  let $V\in C^1_sC^{2}_y (\Bbb R^{3+1})$  be a time periodic solution to $\se$ with the period $S_0 >0$. Suppose there exists an open set $G\subset \Bbb R^3$ containing  the origin such that
  \bb\label{vani}
 V= 0 \quad\mbox{on}\quad G\times (0, S_0).
  \ee

Then,  $V=0$  on $\Bbb R^{3+1}$.
\end{thm}
\noindent{\em Remark 1.1 } We emphasize that our spatial regularity assumption of $V(y,s)$ is only $C^{1,\delta} (\Bbb R^3)$,
and we do not assume any decay condition on $V$ or its derivatives at the spatial infinity.
 Therefore the above result is completely different from \cite{cha4}.\\
\ \\
Let us consider time-periodic solutions of the following more general system than $\se$.
\bb\label{v}
 \left\{\aligned  &\frac{\partial V}{\partial s}  +aV+b (y\cdot\nabla ) V+ (V\cdot \nabla) V =-\nabla P, \\
&\mathrm{div }\, V=0,\\
&V(y,0)=V_0 (y),
\endaligned \right.
\ee
 where $a,b\in \Bbb R$, and $b\neq 0$. For the system (\ref{v}) we shall prove following result in the next section, from which  Theorem 1.1 follows as an immediate corollary. Below we use the notation $B_\rho=\{ y\in \Bbb R^3\,|\, |y|<\rho\}$.
\begin{thm}
Let $a,b \in \Bbb R$ with $b\neq0$, $\delta\in (0, 1)$, and let $V\in C_s ^1C^{1,\delta}_y(\Bbb R^{3+1})$ with $\O=\mathrm{curl}\, V  \in C^1(\Bbb R^{3+1})$ be a time periodic solution to (\ref{v}) with the period $S_0 > 0$. We assume:
 \begin{itemize}
 \item[(i)] Either the case $a+b>0$ and $b>0$, or $a+b<0$ and $b<0$:\\
The origin $y=0$ is a point of local extremum of $V_i (\cdot,s)$ for all $i=1,2,3$ and $s\in (0, S_0)$, and  $V(0,s)=0$.
 \item[(ii)] Either the  case $a+b\leq 0, b>0$, or $a+b\geq 0, b<0$:\\
The origin  $y=0$ is a point of local extremum of $V_i (\cdot,s)$ for all $i=1,2,3$ and $s\in (0, S_0)$, and there exist $\rho>0$ and $\beta > \frac{|a+b|}{b}+2$ such that
\bb\label{regular}
 V\in C_s ^1C_y^{\beta} (B_\rho\times (0,S_0)), \quad \mbox{and}\quad  D^k V(0, s)=0 \quad \forall k=0, \cdots, [\beta],
 \forall s\in (0, S_0).
 \ee
 \end{itemize}
Then, $V=0$ on $\Bbb R^{3+1}$.
\end{thm}
\subsection{The magnetohydrodynamic system}

The result of the previous subsection also holds for more general system such as the following inviscid magnetohydrodynamic equations in $\Bbb R^3$,
\[
\mathrm{ (MHD) }
 \left\{ \aligned
 &\frac{\partial v}{\partial t} +(v\cdot \nabla )v =(b\cdot\nabla)b-\nabla (p +\frac12 |b|^2), \\
 &\frac{\partial b}{\partial t} +(v\cdot \nabla )b =(b \cdot \nabla
 )v,\\
 &\quad \textrm{div }\, v =\textrm{div }\, b= 0 ,\\
  &v(x,0)=v_0 (x), \quad b(x,0)=b_0 (x),
  \endaligned
  \right.
  \]
where $v=(v_1, v_2 , v_3 )$, $v_j =v_j (x, t)$, $j=1,\cdots,3$,
is the velocity of the flow, $p=p(x,t)$ is the scalar pressure,
$b=(b_1, b_2 , b_3 )$, $b_j =b_j (x, t)$, is the magnetic field,
and $v_0$, $b_0$ are the given initial velocity and magnetic field,
 satisfying div $v_0 =\mathrm{div}\, b_0= 0$, respectively. If we set $b=0$ in (MHD), then it reduces to (E).
  The scaling property of (MHD) is the same as (E).
 Given $\alpha \neq -1$ and $T>0$, we  make the  self-similar
transform of (MHD), which is defined by the map $(v,b, p) \mapsto (V,B,P)$
given by
\begin{equation}
\label{31a}
 v(x, t)=\frac {1}{(T-t)^{\frac{\a}{\a+1}}} V(y,s), \,
 b(x, t)=\frac {1}{(T-t)^{\frac{\a}{\a+1}}} B(y,s),\,
 p(x,t)=\frac {1}{(T-t)^{\frac{2\a}{\a+1}}} P(y,s),
\end{equation}
where
\begin{equation}
\label{32a} y = \frac{x}{(T-t)^{\frac{1}{\a+1}}}, \quad s =  \log \left(\frac{T}{T-t}\right).
\end{equation}
Substituting (\ref{31a})--(\ref{32a}) into (MHD), we have the following system in terms of  $(V,B,P)$:
\begin{equation}\label{mhd} \left\{
\aligned
&\frac {\partial V}{\partial s}+ \frac{\a}{\a+1} V +\frac{1}{\a+1}(y \cdot \nabla)V + (V\cdot \nabla )V =
(B\cdot \nabla ) B-\nabla
(P+\frac12 |B|^2),\\
&\frac {\partial B}{\partial s}+ \frac{\a}{\a+1}B +\frac{1}{\a+1}(y \cdot \nabla)B + (V\cdot \nabla )B =(B\cdot \nabla ) V,\\
& \mathrm{div}\,V=0,\,\,\,\mathrm{div}\,B=0, \\
& V(y,0)=V_0 (y)=T^{\frac{\a}{\a+1}}v_0 (T^{\frac{1}{\a+1}}y),\, B(y,0)=B_0 (y)=T^{\frac{\a}{\a+1}}b_0 (T^{\frac{1}{\a+1}}y).
\endaligned \right.
\end{equation}
In \cite{cha4} it is proved that  the time periodic solution of (\ref{mhd}) is trivial under suitable decay conditions of $V$ and $B$ at spatial infinity, which is in itself a substantial improvement of the result in \cite{cha3}.
Below we prove following unique continuation type property for the time periodic solution of (\ref{mhd}), which can be viewed as also a generalization of Theorem 1.1.
\begin{thm}
Let $\a  \neq -1 $, $ \delta \in (0, 1)$ and  let $(V,B)\in C^1_sC^{2}_y (\Bbb R^{3+1})\times C^1(\Bbb R^{3+1})$  be a time periodic solution to (\ref{mhd}) with the period $S_0 > 0$.   We assume the following local vanishing conditions depending on $\a$.
\begin{itemize}
\item[(i)] If $\a >0$ , then we assume
$V(0,s)=0$ on $G \times (0, S_0)$.
\item[(ii)] If $-1<\a \leq 0 $ or $\a <-1$, then we assume there exists an open set $G\subset \Bbb R^3$ containing the origin such that
$V=B=0$ on $G\times (0, S_0)$.
\end{itemize}
Then,  $V=B=0$  on $\Bbb R^{3+1}$.
\end{thm}
Theorem 1.3 is a corollary of the following result with weaker hypothesis.
\begin{thm}
Let $\a  \neq -1 $, $ \delta \in (0, 1)$ and  let $(V,B)\in C^1_sC^{1, \delta}_y (\Bbb R^{3+1})\times C^1(\Bbb R^{3+1})$ with $\O=\mathrm{curl}\, V  \in C^1(\Bbb R^{3+1})$  be a time periodic solution to (\ref{mhd}) with the period $S_0 > 0$.  We assume the following at the point $y=0$.
\begin{itemize}
\item[(i)] If $\a >0$, we assume that $V(0,s)=0$, and $y=0$ is a point of local extremum of $V_i (\cdot,s)$ for all $i=1,2,3$ and $s\in (0, S_0)$.
\item[(ii)] If $-1<\a \leq 0 $ or $\a <-1$,  we assume that $y=0$ is a point of local extremum of $V_i (\cdot,s)$ for all $i=1,2,3$ and $s\in (0, S_0)$. Moreover, there exists $\rho>0$ and $\beta > |\a|+2$ such that
\bb\label{regulara}
 V\in C_s ^1C_y^{\beta} (B_\rho\times (0,S_0)), \quad \mbox{and}\quad  D^k V(0, s)=0 \quad \forall k=0, \cdots, [\beta],
 \forall s\in (0, S_0).
 \ee
\end{itemize}
Then,  $V=B=0$  on $\Bbb R^{3+1}$.
\end{thm}
\section{ Proof of the Main Theorems }
\setcounter{equation}{0}

\noindent{\bf Proof of Theorem 1.2: }\\
\noindent{\em \underline{(A) The case  with $a+b>0, b>0$:} }
We choose  $r_0>0$  such that $y=0$ is an extremal point of $V_i (y,s)$ on $B_{r_0}$ for each $i=1,2,3$ and $s\in (0, S_0)$.
Taking curl on (\ref{v}), we obtain the vorticity equations,
\bb\label{vol}
\left\{\aligned &\frac {\partial \O}{\partial s}+ (a+b)\O + b (y\cdot \nabla) \O+(V\cdot \nabla )\O =(\O\cdot \nabla )V,\\
& \mathrm{curl} \,V=\O, \quad \mathrm{div} \, V=0,\\
& \O (y,s)=\O_0 (y).\endaligned \right.
\ee
For each $\rho>0$  we  define a bump  function $\psi=\psi_\rho(y)$  as follows.
\bb\label{cut}
\psi_\rho (y)=\left\{ \aligned &0,  &\mbox{if}\quad |y|\geq \rho\\
&\exp\left(\frac{1}{|y|^2-\rho^2}\right), &\mbox{if}\quad |y|<\rho .
                                \endaligned
                                \right.
    \ee
For   $f\in C^\delta (\Bbb R^3)$ and  $\mathcal{D }\subset \Bbb R^3$, a bounded set, we denote
    $$ \|f\|_{C^{\delta} (\mathcal{D})} =\sup_{y_1, y_2\in\mathcal{ D}, y_1\neq y_2}\frac{|f(y_1)-f(y_2)|}{|y_1-y_2|^\delta}.
    $$
Let $R\in (r_0, \infty)$ be any sufficiently large number. We define
 \bb\label{defm}
 M=\sup_{s\in (0, S_0)} \|D V(\cdot,s)\|_{C^{\delta} (B_R)},
 \ee
 and set
 \bb\label{rad1}
 R_0=\min\{r_0, \mu^{\frac{1}{\delta}}\}, \qquad \mu:=\min\left\{ \frac{ b}{2M}, \frac{a+b}{2M}\right\}.
 \ee
Let $y\in B_{R_0}$.  Then, since $|V(0,s)|=|DV(0,s)|=0$ for $s\in (0,S_0)$, we have
\bq\label{nvest}
 |D V (y,s)|&=& |DV(y,s)- DV(0 ,s)|\leq \sup_{s\in (0, S_0)}\|DV(\cdot, s)\|_{C^{\delta} (B_R )} |y|^\delta\n \\
&\leq& R_0^\delta M=\mu M\leq \frac{a+b}{2},
\eq
and
\bq\label{vest}
|V(y,s)|&=&\left| V(0, s)+\int_0 ^1 y\cdot \nabla V(\tau y,s)d\tau\right|\n \\
 &\leq&  |y| \sup_{|z|\leq |y|} |D V(z,s)|\leq |y|^{1+\delta}\sup_{s\in (0, S_0)}\|DV(\cdot,s)\|_{C^{\delta} (B_R )}  \n \\
  &\leq&   |y|R_0^\delta  M\leq |y|\mu M\leq \frac{b}{2} |y|
\eq
for all $(y,s)\in B_{R_0}\times (0, S_0)$.  Multiplying  (\ref{vol}) by $\O \psi_{R_0} $, we obtain
\bq\label{21}
\lefteqn{\frac{\partial}{\partial s} (\psi_{R_0}|\O|^2 ) + 2(a+b)\psi_{R_0} |\O|^2 +b(y\cdot \nabla )( \psi_{R_0}|\O|^2 ) +(V\cdot \nabla )(\psi_{R_0}|\O|^2 )}\hspace{.0in}\n \\
 &&\leq 2 \psi_{R_0}|\nabla V ||\O|^2+ b |\O|^2 (y\cdot \nabla )\psi_{R_0} +|\O|^2 (V\cdot \nabla )\psi_{R_0}.
 \eq

Since $\psi_{R_0}(\cdot )$ is radially non-increasing, we have from (\ref{vest}) and (\ref{nvest}) respectively
\bb\label{25}
 |(V\cdot \nabla )\psi_{R_0}|\leq -|V| \frac{\partial \psi_{R_0}}{\partial r}\leq -\frac{b}{2}|y| \frac{\partial \psi_{R_0}}{\partial r}= -\frac{b}{2}   (y\cdot \nabla) \psi_{R_0},
\ee
and
 \bb\label{24}
|\nabla V(y,s)|\psi_{R_0}
\leq \frac{(a+b)}{2} \psi_{R_0}(y),
 \ee
  for all $(y,s)\in B_{R_0}\times (0, S_0)$.
Substituting (\ref{24}) and (\ref{25})  into (\ref{21}),  we  have
 \bq\label{26}
\lefteqn{\frac{\partial}{\partial s} (|\O|^2\psi_{R_0}) +(a+b) |\O|^2\psi_{R_0} +b(y\cdot \nabla )(|\O|^2 \psi_{R_0}) }\hspace{1.5in}\n \\
&&+(V\cdot \nabla )(|\O|^2\psi_{R_0})
\leq \frac{b}{2}  |\O|^2(y\cdot \nabla )\psi_{R_0}
 \eq
 for all  $(y,s)\in B_{R_0}\times (0, S_0)$. Since  $\frac{ b}{2} |\O|^2(y\cdot \nabla )\psi_{R_0} \leq 0$, we  obtain from (\ref{26}) the following differential inequality.
 \bb\label{27}
 \frac{\partial}{\partial s} f(y,s)+(a+b) f(y,s)+b y\cdot \nabla f(y,s)+V\cdot \nabla f(y,s)\leq 0,
 \ee
 where we set
 $
 f(y,s):= |\O|^2\psi_{R_0}.
 $
 Let $q>2$. Multiplying (\ref{27}) by $f|f|^{q-2}$, and integrate it over $B_{R_0}\times (0, S_0)$, and taking into account the fact $f(y,0)=f(y, S_0)$ for all $y\in B_{R_0}$, one has
 $$
 \left(a+b-\frac{3b}{q} \right) \int_0 ^{S_0} \int_{B_{R_0}} |f|^q dyds \leq 0.
 $$
 Choosing $q>\frac{3b}{a+b}$, we have $f=0$, and therefore $\O=0$ on $B_{R_0}\times (0, S_0)$.
 Since $\mathrm{div} \, V=0$ and $\mathrm{curl}\, V=0$, we find that each  $V_i (\cdot,s)$, $i=1,2,3$, is a harmonic function on  $B_{R_0}$. Since $y=0$ is an extremal point of $V_i (\cdot,s)$ on $B_{R_0}$, and $V_i(0,s)=0$ for each $s\in (0, S_0)$, we have $V_i=0$  on $B_{R_0}\times (0,S_0)$  for each $i=1,2,3$  by the maximum principle for  harmonic functions.\\

 Now we will inductively extend the zero set of $V$ from $B_{R_0}$ to a sequence of increasing balls.
Let $m$ be the smallest integer such that
 $ \frac{R-R_0}{\mu^{\frac{1}{\delta }}}\leq m$.
 Consider the sequence $\{R_k\}_{k=1}^{m}$, where
 $$R_{k}= R_{k-1} +\mu^{\frac{1}{\delta}} , \qquad k=1, \cdots, m.$$
Given $y\in B_{R_k}$, we define  $\bar{y}=\frac{y}{|y|} R_{k-1}\in \partial B_{R_{k-1}}$ for $k=1, \cdots, m$. Suppose we have shown that
 $V=0$ on $B_{R_{k-1}} \times (0, S_0)$, which is the case for $k=1$.
 Then, instead of (\ref{nvest}) and (\ref{vest}) we have
 \bq\label{nvesta}
|\nabla V(y,s)|&\leq& |y-\bar{y}|^\delta \sup_{s\in (0, S_0)}\|D V(\cdot,s)\|_{C^{\delta} (B_R)}\n \\
&\leq& (R_k-R_{k-1})^\delta M\leq \mu M\leq \frac{a+b}{2},
\eq
and
\bq\label{vesta}
|V(y,s)|
 &\leq&  |y-\overline{y}| \sup_{|z|\leq |y|} |\nabla V(z,s)|\n \\
  &\leq& |y-\overline{y}|^{1+\delta} M \leq  |y-\overline{y}|(R_k-R_{k-1})^\delta M\n \\
  &\leq& |y|\mu M\leq \frac{b}{2} |y|
\eq
for all $(y,s)\in B_{R_k}\times (0, S_0)$.  Multiplying  (\ref{vol}) by $\O \psi_{R_k} $, we obtain
\bq\label{233}
\lefteqn{\frac{\partial}{\partial s} (\psi_{R_k}|\O|^2 ) + 2(a+b)\psi_{R_k} |\O|^2 +b(y\cdot \nabla )( \psi_{R_k}|\O|^2 ) +(V\cdot \nabla )(\psi_{R_k}|\O|^2 )}\hspace{.0in}\n \\
 &&\leq 2 \psi_{R_k}|\nabla V ||\O|^2+ b |\O|^2 (y\cdot \nabla )\psi_{R_k} +|\O|^2 (V\cdot \nabla )\psi_{R_k}.
 \eq

We have from (\ref{vesta}) and (\ref{nvesta}) respectively
\bb\label{255}
 |(V\cdot \nabla )\psi_{R_k}|\leq -|V| \frac{\partial \psi_{R_k}}{\partial r}\leq -\frac{b}{2}|y| \frac{\partial \psi_{R_k}}{\partial r}= -\frac{b}{2}   (y\cdot \nabla) \psi_{R_k},
\ee
and
 \bb\label{244}
|\nabla V(y,s)|\psi_{R_k}
\leq \frac{(a+b)}{2} \psi_{R_k}(y)
 \ee
  for all $(y,s)\in B_{R_k}\times (0, S_0)$.
Substituting (\ref{244}) and (\ref{255})  into (\ref{233}),  we  have
 \bq\label{266}
\lefteqn{\frac{\partial}{\partial s} (|\O|^2\psi_{R_k}) +(a+b) |\O|^2\psi_{R_k} +b(y\cdot \nabla )(|\O|^2 \psi_{R_k}) }\hspace{1.5in}\n \\
&&+(V\cdot \nabla )(|\O|^2\psi_{R_k})\leq \frac{ b}{2} |\O|^2(y\cdot \nabla )\psi_{R_k} \leq 0
 \eq
 for all  $(y,s)\in B_{R_k}\times (0, S_0)$.
 Let $
 f(y,s):= |\O|^2\psi_{R_k},
 $
 and $q>2$. Multiplying (\ref{27}) by $f|f|^{q-2}$, and integrate it over $B_{R_k}\times (0, S_0)$, and taking into account the fact $f(y,0)=f(y, S_0)$ for all $y\in B_{R_k}$, one has
 $$
 \left(a+b-\frac{3b}{q} \right) \int_0 ^{S_0} \int_{B_{R_k}} |f|^q dyds \leq 0.
 $$
 Choosing $q>\frac{3b}{a+b}$, we have $f=0$, and therefore $\O=0$ on $B_{R_k}\times (0, S_0)$. Combining this with the condition $\mathrm{div }\, V=0$, we have that $V(\cdot, s)$ is harmonic, and thus real analytic on $B_{R_k}$, which is zero on $B_{R_{k-1}}$.
 Thus we can extend the zero set of $V(\cdot, s)$ from $B_{R_{k-1}}$ to $B_{R_k}$ for all $s\in (0, S_0)$.
 One can do the above argument for $k=1, \cdots, m$ to show that   $V=0$ on $B_{R-\mu^\frac{1}{\delta}}\times (0, S_0)$. Since $R\in (R_0, \infty)$ can be arbitrarily large, we have shown $V=0$
 on $\Bbb R^3 \times (0, S_0)$.\\
 \ \\
\noindent{\em \underline{(B) The case with $a+b<0, b<0$:} }
 In this case we define $\overline{V} (y, s)= V(y, S_0-s)$,  $\bar{P} (y, s)= P(y, S_0-s)$ and $\Ob (y, s)= \O(y, S_0-s)$ for $0\leq s\leq S_0$. Then, the vorticity equation becomes
 \bb\label{211}
\frac {\partial \Ob}{\partial s}-(a+b)\Ob - b (y\cdot \nabla) \Ob-(\vb\cdot \nabla )\Ob =-(\overline{\O}\cdot \nabla )\vb.
\ee
 This is the same situation as (A) above, since $-(a+b)>0, -b>0$. In particular, we note that the signs in front of the terms
 $(V\cdot \nabla )\O$ and $ (\O\cdot \nabla )V$ are not important in the estimates (\ref{24}) and (\ref{25}).
 Repeating the argument of (A) word by word,  we conclude $V=0$ on $\Bbb R^3\times(0, S_0)$. \\
  \ \\
 \noindent{\em \underline{(C) The case  with $a+b\leq 0, b>0 $:} }
 We choose  $r_0>0$  such that $y=0$ is an extremal point of $V_i (y,s)$ on $B_{r_0}$ for each $i=1,2,3$ and $s\in (0, S_0)$.
Let $R\in (r_0, \infty)$ be fixed.  Similarly to the  proof (A) we define
$M$ as in  (\ref{defm}).
We choose $\gamma < \frac{a+b}{b}$.
Then, we set
\bb\label{mu1}
R_0= \min\{r_0, \mu^{\frac{1}{\delta}}\},  \qquad \mu:=\min\left\{\frac{b}{2 M}, \frac{a+b-b\gamma}{4M}, \frac{a+b-b\gamma}{4b|\gamma| M}\right\}.
\ee
 We  multiply  (\ref{vol}) by $\O \psi_{R_0} |y|^{2\gamma} $ to obtain
\bq\label{212a}
\lefteqn{\frac{\partial}{\partial s} (\psi_{R_0}|y|^{2\gamma}|\O|^2 ) + 2(a+b-b\gamma)\psi_{R_0}|y|^{2\gamma} |\O|^2 +b(y\cdot \nabla )( \psi_{R_0}|y|^{2\gamma}|\O|^2 ) }\hspace{.5in}\n \\
 &&\qquad+(V\cdot \nabla )(\psi_{R_0}|y|^{2\gamma}|\O|^2 )\leq 2|y|^{2\gamma}\psi_{R_0} |\nabla V| |\O|^2+ b |y|^{2\gamma}|\O|^2 (y\cdot \nabla )\psi_{R_0} \n \\
 &&\qquad\quad\qquad+|y|^{2\gamma}|\O|^2 (V\cdot\nabla)\psi_{R_0} +2b\gamma|y|^{2\gamma-1} |\O|^2 V\cdot\frac{y}{|y|}\psi_{R_0}.
 \eq
 We note that thanks to the assumption (\ref{regular}) each terms involving $|y|^{2\gamma} |\O|^2$ belong to $C_y ^1 (B_{R_0}\times (0, S_0))$.
We estimate
 from (\ref{nvest}) and (\ref{vest}),
 \bb\label{213a}
|\nabla V(y,s)|\psi_{R_0} \leq  |y|^\delta M \psi_{R_0}\leq R_0^\delta M \psi_{R_0} \leq \frac14 (a+b-b\gamma)\psi_{R_0},
 \ee
\bq\label{214a}
 |(V\cdot \nabla )\psi_{R_0}|&\leq& - |V|\partial_r \psi_{R_0} \leq  -M|y|^{1+\delta} \partial_r \psi_{R_0} \leq -M R_0^\delta |y|\partial_r \psi_{R_0} \n \\
 &\leq& -\frac{b}{2}  |y| \partial_r \psi_{R_0}
 =  -\frac{b}{2} (y\cdot \nabla) \psi_{R_0},
\eq
and
\bq\label{214aa}
2|\gamma|b |y|^{2\gamma-1}  \left|V\cdot \frac{y}{|y|} \right|\psi_{R_0}&\leq& 2|\gamma|b M|y|^{2\gamma-1}|y|^{1+\delta} \psi_{R_0} \leq 2b|\gamma|R_0^\delta M|y|^{2\gamma} \psi_{R_0} \n \\
&\leq&
\frac12 (a+b-b\gamma) |y|^{2\gamma} \psi_{R_0}
\eq
for all $s\in (0, S_0)$.
Substituting (\ref{213a})-(\ref{214aa})  into (\ref{212a}),  we  have
 \bq\label{215a}
\lefteqn{\frac{\partial}{\partial s} (\psi_{R_0} |y|^{2\gamma} |\O|^2 ) +2(a+b-b\gamma) \psi_{R_0} |y|^{2\gamma} |\O|^2  +b(y\cdot \nabla )(\psi_{R_0} |y|^{2\gamma} |\O|^2 ) +(V\cdot \nabla )(\psi_{R_0} |y|^{2\gamma}|\O|^2 )}\hspace{.1in}\n \\
&&\leq 2 \psi_{R_0} |y|^{2\gamma}|\nabla V(y,s)| |\O|^2 +b |y|^{2\gamma}|\O|^2 (y\cdot \nabla )\psi_{R_0}\n \\
 &&\hspace{1.in}\qquad+|y|^{2\gamma}|\O|^2 (V\cdot \nabla )\psi_{R_0}  +2\gamma b|y|^{2\gamma-1} |\O|^2 V\cdot \frac{y}{|y|} \psi_{R_0} \n \\
&& \leq \frac12 (a+b-b\gamma) \psi_{R_0} |y|^{2\gamma}|\O|^2+ \frac{b}{2}  |\O|^2|y|^{2\gamma}(y\cdot \nabla )\psi_{R_0}
+\frac12 (a+b-b\gamma)\psi_{R_0} |y|^{2\gamma} |\O|^2\n \\
 &&=(a+b-b\gamma) \psi_{R_0} |y|^{2\gamma}|\O|^2+\frac{b}{2}  |\O|^2|y|^{2\gamma}(y\cdot \nabla ) \psi_{R_0}
  \eq
 for all  $s\in (0, S_0)$. Since $\frac{ b}{2}|y|^{2\gamma} |\O|^2(y\cdot \nabla )\psi_{R_0}  \leq 0$, we have the following differential inequality from (\ref{215a}),
 \bb\label{216a}
 \frac{\partial}{\partial s} f(y,s)+(a+b-b\gamma) f(y,s)+b y\cdot \nabla f(y,s)+V\cdot \nabla f(y,s)\leq 0,
 \ee
 where we set
 $
 f(y,s):=\psi_{R_0} |y|^{2\gamma} |\O|^2.
 $
 For $q>2$ we multiply (\ref{216}) by $f|f|^{q-2}$, and integrate it over $B_{R_0}\times (0, S_0)$, and using the fact that $f(y,0)=f(y,S_0)$ for all $y\in B_{R_0}$, we obtain
 $$
 \left(a+b-b\gamma-\frac{3b}{q}\right ) \int_0 ^{S_0} \int_{B_{R_k}} |f(y,s)|^q dyds \leq 0.
 $$
 Choosing $q> \frac{3b}{a+b-b\gamma}$, we have $f=0$, and therefore $\O=0$ on $B_{R_0}\times (0, S_0)$. Thus, $V(\cdot, s)$ is harmonic on  $B_{R_0}$. Since $y=0$ is an extremal point of $V_i (\cdot,s)$ on $B_{R_0}$, and $V_i(0,s)=0$ for each $s\in (0, S_0)$, we have $V_i=0$  on $B_{R_0}\times (0,S_0)$  for each $i=1,2,3$  by the maximum principle for the harmonic function.\\

Next we extend the zero set of $V$ from $B_{R_0}\times (0, S_0)$ to successively increasing sets by induction argument.
We define
\bb\label{rad2}
R_{k} =R_{k-1}+\mu^{\frac{1}{\delta}}, \qquad k=1, \cdots, m
 \ee
 where $\mu$ is the same number as in (\ref{mu1}), and $m$ is the smallest integer such that $\frac{R-R_0}{\mu^{\frac{1}{\delta}}} \leq m$.
 We  multiply  (\ref{vol}) by $\O \psi_{R_k} |y|^{2\gamma} $ to obtain
\bq\label{212}
\lefteqn{\frac{\partial}{\partial s} (\psi_{R_k}|y|^{2\gamma}|\O|^2 ) + 2(a+b-b\gamma)\psi_{R_k}|y|^{2\gamma} |\O|^2 +b(y\cdot \nabla )( \psi_{R_k}|y|^{2\gamma}|\O|^2 ) }\hspace{.5in}\n \\
 &&\qquad+(V\cdot \nabla )(\psi_{R_k}|y|^{2\gamma}|\O|^2 )\leq 2|y|^{2\gamma}\psi_{R_k} |\nabla V| |\O|^2+ b |y|^{2\gamma}|\O|^2 (y\cdot \nabla )\psi_{R_k} \n \\
 &&\qquad\quad\qquad+|y|^{2\gamma}|\O|^2 (V\cdot\nabla)\psi_{R_k} +2b\gamma|y|^{2\gamma-1} |\O|^2 V\cdot\frac{y}{|y|}\psi_{R_k}.\n \\
 \eq
 We suppose that the proof of $V=0$ on $B_{k-1}\times (0, S_0)$ is done, which is the case for $k=1$.

For $y\in  B_{R_k}$ let $\overline{y}=\frac{y}{|y|} R_{k-1}\in \partial B_{R_{k-1}}$.
 We estimate
 from (\ref{nvest}) and (\ref{vest}),
 \bb\label{213}
|\nabla V(y,s)|\psi_{R_k} \leq  |y-\bar{y}|^\delta M \psi_{R_k}\leq (R_k-R_{k-1})^\delta M \psi_{R_k} \leq \frac14 (a+b-b\gamma)\psi_{R_k},
 \ee
\bq\label{214}
 |(V\cdot \nabla )\psi_{R_k}|&\leq& - |V|\partial_r \psi_{R_k} \leq  -M|y-\bar{y}|^{1+\delta} \partial_r \psi_{R_k} \leq -M(R_k-R_{k-1})^\delta |y|\partial_r \psi_{R_k} \n \\
 &\leq& -\frac{b}{2}  |y| \partial_r \psi_{R_k}
 =  -\frac{b}{2} (y\cdot \nabla) \psi_{R_k},
\eq
and
\bq\label{214a}
2|\gamma|b |y|^{2\gamma-1}  \left|V\cdot \frac{y}{|y|} \right|\psi_{R_k}&\leq& 2|\gamma|b M|y|^{2\gamma-1}|y-\bar{y}|^{1+\delta} \psi_{R_k} \leq 2b|\gamma|(R_k-R_{k-1})^\delta M|y|^{2\gamma} \psi_{R_k} \n \\
&\leq&
\frac12 (a+b-b\gamma) |y|^{2\gamma} \psi_{R_k}
\eq
for all $s\in (0, S_0)$.
Substituting (\ref{213})-(\ref{214a})  into (\ref{212}),  we  have
 \bq\label{215}
\lefteqn{\frac{\partial}{\partial s} (\psi_{R_k} |y|^{2\gamma} |\O|^2 ) +2(a+b-b\gamma) \psi_{R_k} |y|^{2\gamma} |\O|^2  +b(y\cdot \nabla )(\psi_{R_k} |y|^{2\gamma} |\O|^2 ) +(V\cdot \nabla )(\psi_{R_k} |y|^{2\gamma}|\O|^2 )}\hspace{.1in}\n \\
&&\leq 2 \psi_{R_k} |y|^{2\gamma}|\nabla V(y,s)| |\O|^2 +b |y|^{2\gamma}|\O|^2 (y\cdot \nabla )\psi_{R_k}\n \\
 &&\hspace{1.in}\qquad+|y|^{2\gamma}|\O|^2 (V\cdot \nabla )\psi_{R_k}  +2\gamma b|y|^{2\gamma-1} |\O|^2 V\cdot \frac{y}{|y|} \psi_{R_k} \n \\
&& \leq \frac12 (a+b-b\gamma) \psi_{R_k} |y|^{2\gamma}|\O|^2+ \frac{b}{2}  |\O|^2|y|^{2\gamma}(y\cdot \nabla )\psi_{R_k}
+\frac12 (a+b-b\gamma)\psi_{R_k} |y|^{2\gamma} |\O|^2\n \\
 &&=(a+b-b\gamma) \psi_{R_k} |y|^{2\gamma}|\O|^2+\frac{b}{2}  |\O|^2|y|^{2\gamma}(y\cdot \nabla ) \psi_{R_k}
  \eq
 for all  $s\in (0, S_0)$. Since $\frac{ b}{2}|y|^{2\gamma} |\O|^2(y\cdot \nabla )\psi_{R_k}  \leq 0$, we have the following differential inequality from (\ref{215}),
 \bb\label{216}
 \frac{\partial}{\partial s} f(y,s)+(a+b-b\gamma) f(y,s)+b y\cdot \nabla f(y,s)+V\cdot \nabla f(y,s)\leq 0,
 \ee
 where we set
 $
 f(y,s):=\psi_{R_k} |y|^{2\gamma} |\O|^2.
 $
 For $q>2$ we multiply (\ref{216}) by $f|f|^{q-2}$, and integrate it over $B_{R_k}\times (0, S_0)$, and using the fact that $f(y,0)=f(y,S_0)$ for all $y\in B_{R_k}$, we obtain
 $$
 \left(a+b-b\gamma-\frac{3b}{q}\right ) \int_0 ^{S_0} \int_{B_{R_k}} |f(y,s)|^q dyds \leq 0.
 $$
 Choosing $q> \frac{3b}{a+b-b\gamma}$, we have $f=0$, and therefore $\O=0$ on $B_{R_k}\times (0, S_0)$. Thus, $V(\cdot, s)$ is harmonic, and hence real analytic on $B_{R_k}$, and zero on $B_{R_{k-1}}$ for each $s\in (0, S_0)$. By extension of the zero set for real analytic function we obtain $V=0$ on $B_{R_k}\times (0, S_0)$.
 Similarly to the proof (A), we can repeat the above argument  for $k=1, \cdots, m$, and we reach the domain $B_{R-\mu^\frac{1}{\delta}}\times (0, S_0)$ for the zero set of $V$. Since $R$ is arbitrary, we are done.\\
 \ \\
\noindent{\em \underline{(D) The case of (ii) with $a+b\geq 0, b<0 $:} }
Similarly to the proof (B) above  we introduce  $\overline{V} (y, s)= V(y, S_0-s)$,  $\bar{P} (y, s)= P(y, S_0-s)$ and $\Ob (y, s)= \O(y, S_0-s)$ for $0\leq s\leq S_0$ to derive (\ref{211}).  Then we are reduced to the case of (C) by similar argument to (B). $\square$\\

\noindent{\bf Proof of Theorem 1.4 :  } \\
\noindent{\em \underline{(A) The case  with $\a >0$:} }
We choose  $r_0>0$  such that $y=0$ is an extremal point of $V_i (y,s)$ on $B_{r_0}$ for each $i=1,2,3$ and $s\in (0, S_0)$.
Let $R\in (r_0, \infty)$ be any sufficiently large number.
Let us set
 \bb\label{rad3}
R_0=\min\{ r_0, \mu ^{\frac{1}{\delta}}\}, \qquad \mu=\min\left\{\frac{1}{2(\a+1) M},\, \frac{\a}{2(\a+1) M}\right\},
 \ee
 where $M$ is the same as in (\ref{defm}).
Let $y\in B_{R_0}$.  Then, since $|V(0,s)|=|DV(0,s)|=0$ for $s\in (0,S_0)$, we have
\bq\label{nvestb}
 |D V (y,s)|&=& |DV(y,s)- DV(0 ,s)|\leq \sup_{s\in (0, S_0)}\|DV(\cdot, s)\|_{C^{\delta} (B_R )} |y|^\delta\n \\
&\leq& R_0^\delta M=\mu M\leq \frac{\a}{2(\a+1)},
\eq
and
\bq\label{vestb}
|V(y,s)|&=&\left| V(0, s)+\int_0 ^1 y\cdot \nabla V(\tau y,s)d\tau\right|\n \\
 &\leq&  |y| \sup_{|z|\leq |y|} |D V(z,s)|\leq |y|^{1+\delta}\sup_{s\in (0, S_0)}\|DV(\cdot,s)\|_{C^{\delta} (B_R )}  \n \\
  &\leq&   |y|R_0^\delta  M\leq |y|\mu M\leq \frac{1}{2(\a+1)} |y|
\eq
for all $(y,s)\in B_{R_0}\times (0, S_0)$.
 Multiplying  the second equations of (\ref{mhd}) by $B \psi_{R_0} $, we obtain
\bq\label{30b}
\lefteqn{\frac{\partial}{\partial s} (\psi_{R_0}|B|^2 ) + \frac{2\a}{\a+1}\psi_{R_0} |B|^2 +\frac{1}{\a+1}(y\cdot \nabla )( \psi_{R_0}|B|^2 ) +(V\cdot \nabla )(\psi_{R_0}|B|^2 )}\hspace{.0in}\n \\
 &&\leq 2 \psi_{R_0}|\nabla V ||B|^2+\frac{1}{\a+1} |B|^2 (y\cdot \nabla )\psi_{R_0} +|B|^2 (V\cdot \nabla )\psi_{R_0}.
 \eq
 We have from (\ref{vestb})
\bb\label{25bb}
 |(V\cdot \nabla )\psi_{R_0}|\leq -|V| \frac{\partial \psi_{R_0}}{\partial r}\leq -\frac{1}{2(\a+1)}|y| \frac{\partial \psi_{R_0}}{\partial r}= -\frac{1}{2(\a+1)}   (y\cdot \nabla) \psi_{R_0},
\ee

Substituting (\ref{nvestb}) and (\ref{25bb}) into (\ref{30b}), we have
\bq\label{30b}
\lefteqn{\frac{\partial}{\partial s} (\psi_{R_0}|B|^2 ) + \frac{\a}{\a+1}\psi_{R_0} |B|^2 +\frac{1}{\a+1}(y\cdot \nabla )( \psi_{R_0}|B|^2 ) +(V\cdot \nabla )(\psi_{R_0}|B|^2 )}\hspace{2.in}\n \\
 &&\leq \frac{1}{2(\a+1)} |B|^2 (y\cdot \nabla )\psi_{R_0}\leq 0
 \eq
Let $q>2$. Let $f=\psi_{R_0}|B|^2$.  Multiplying  (\ref{30b}) by $f |f|^{q-2}$, and integrating  it over $B_{R_0}\times (0, S_0)$,
we find that
$$ \frac{1}{\a+1} \left(\a-\frac{3}{q}\right)\int_0 ^{S_0} \int_{B_{R_0}}| f|^qdyds \leq 0. $$
Choosing $q>3/\a$, we obtain $f=0$, and therefore $B=0$ on $B_{R_0}\times (0, S_0)$. Hence, the first equation of (\ref{mhd}) reduces to $\se$  on $B_{R_0}\times (0, S_0)$.  By repeating the proof of Theorem 1.2 (A) for the system $\se$, we conclude
$V=0$ also on $B_{R_0}\times (0, S_0)$.\\

Let $R\in (R_0, \infty)$ be any sufficiently large number. As previously we will use induction argument to show that
$V=B=0$ on $B_{R-\mu^{\frac{1}{\delta}}}\times (0, S_0)$.
Let $m$ is the smallest integer such that $\frac{R-R_0}{\mu^{\frac{1}{\delta}}} \leq m$, and we define the sequence $\{ R_k\}_{k=1}^{m}$ by $R_{k}=R_{k-1} +\mu^{\frac{1}{\delta}}, k=1, \cdots ,m$.
We denote $\bar{y}=\frac{y}{|y|}R_{k-1}$ for $y\in B_{R_k}$.
We suppose $V=B=0$ on $B_{R_{k-1}}\times (0, S_0)$, which is the case for $k=1$.
Then, we multiply  the second equations of (\ref{mhd}) by $B \psi_{R_k}$ to obtain
\bq\label{31}
\lefteqn{\frac{\partial}{\partial s} (\psi_{R_k}|B|^2 ) + \frac{2\a}{\a+1}\psi_{R_k} |B|^2 +\frac{1}{\a+1}(y\cdot \nabla )( \psi_{R_k}|B|^2 ) +(V\cdot \nabla )(\psi_{R_k}|B|^2 )}\hspace{.0in}\n \\
 &&\leq 2 \psi_{R_k}|\nabla V ||B|^2+ \frac{1}{\a+1} |B|^2 (y\cdot \nabla )\psi_{R_k} +|B|^2 (V\cdot \nabla )\psi_{R_k}.
 \eq
From (\ref{vest}) and (\ref{nvest}) we estimate
\bq\label{35}
 |(V\cdot \nabla ) \psi_{R_k}|&\leq& -|V|  \frac{\partial  \psi_{R_k}}{\partial r} \leq -M |y-\bar{y}|^{1+\delta} \frac{\partial  \psi_{R_k}}{\partial r}\n \\
 &\leq& -M (R_k-R_{k-1})^\delta|y| \frac{\partial \psi_{R_k}}{\partial r}\leq  -\frac{1}{2(\a+1)}   (y\cdot \nabla) \psi_{R_k},
\eq
and
 \bb\label{34}
|\nabla V(y,s)| \psi_{R_k}\leq |y-\bar{y}|^\delta M  \psi_{R_k}  \leq (R_k-R_{k-1})^\delta M \psi_{R_k}
\leq \frac{\a}{2(\a+1)}  \psi_{R_k},
 \ee
  for all $s\in (0, S_0)$.
Substituting (\ref{35}) and (\ref{34})  into (\ref{31}),  we  have
 \bq\label{36}
\lefteqn{\frac{\partial}{\partial s} (|B|^2 \psi_{R_k}) +\frac{\a}{\a+1} |B|^2 \psi_{R_k} +\frac{1}{\a+1}(y\cdot \nabla )(|B|^2  \psi_{R_k}) }\hspace{1.5in}\n \\
&&+(V\cdot \nabla )(|B|^2 \psi_{R_k})
\leq \frac{b}{2}  |B|^2(y\cdot \nabla ) \psi_{R_k}\leq 0,
 \eq
 and
 \bb\label{36a}
\frac{ \partial f}{\partial s} +\frac{\a}{\a+1}f +\frac{1}{\a+1}(y\cdot \nabla )f+V\cdot \nabla f\leq 0, \qquad f:=|B|^2 \psi_{R_k}.
\ee
 for all  $(y,s)\in B_{R_k}\times(0, S_0)$.
 Multiplying (\ref{36a}) by $f |f|^{q-2}$, and integrating it over $B_{R_k}\times (0, S_0)$, and choosing $q>3/\a$, we have $f=0$, and $B=0$ on $B_{R_k}\times (0, S_0)$.
 Therefore, the system (\ref{mhd}) reduces to $\se$ on $B_{R_k}\times (0, S_0)$, and $V=0$ on $B_{R_{k-1}}\times (0, S_0)$.
 By Theorem 1.2 and, in particular, by the proof (A) we can extend the zero set of $V$, which is a solution to
  $\se$, up to $B_{R_k}\times(0, S_0)$.  Therefore we have extended the zero set of $V$ and $B$  for (\ref{mhd})  up to
  $ B_{R_k}\times(0, S_0)$. Then, we repeat the above argument for $k=1, \cdots, m$ to have the zero set of both $V$ and $B$ as $B_{R-\mu^\frac{1}{\delta}}\times (0, S_0)$, and since $R\in (R_0, \infty)$ is arbitrary, the proof is complete.\\

\noindent{\em \underline{(B) The case  with $-1< \a \leq 0 $:} }
We choose  $r_0>0$  such that $y=0$ is an extremal point of $V_i (y,s)$ on $B_{r_0}$ for each $i=1,2,3$ and $s\in (0, S_0)$.
Let $R\in (r_0, \infty)$ be any sufficiently large number. We choose $\gamma < \a $, and
 set
 \bb\label{rad3}
R_0=\min\{ r_0, \mu ^{\frac{1}{\delta}}\}, \qquad \mu:= \min\left\{\frac{1}{2(\a+1) M}, \frac{\a-\gamma}{4|\gamma| M}, \frac{\a-\gamma}{4(\a+1)M}\right\},
 \ee
 where $M$ is the same as in (\ref{defm}).
Let $y\in B_{R_0}$.
 We  multiply  the second equations of (\ref{mhd}) by $B\psi_{R_0}|y|^{2\gamma} $ to obtain
\bq\label{312bb}
\lefteqn{\frac{\partial}{\partial s} (\psi_{R_0}|y|^{2\gamma}|B|^2 ) + \frac{2(\a-\gamma)}{\a+1}\psi_{R_0}|y|^{2\gamma} |B|^2 +\frac{1}{\a+1}(y\cdot \nabla )( \psi_{R_0}|y|^{2\gamma}|B|^2 ) )}\n \\
 &&+(V\cdot \nabla )(\psi_{R_0}|y|^{2\gamma}|B|^2) \leq 2|y|^{2\gamma}\psi_{R_0} |\nabla V| |B|^2+ \frac{1}{\a+1} |y|^{2\gamma}|B|^2 (y\cdot \nabla )\psi_{R_0} \n \\
 &&\hspace{0.in}\qquad\quad\qquad+|y|^{2\gamma}|B|^2 (V\cdot\nabla) \psi_{R_0} +\frac{2\gamma}{\a+1} |y|^{2\gamma-1} |B|^2 V\cdot\frac{y}{|y|}\psi_{R_0}.
 \eq
 We note that thanks to the hypothesis (\ref{regulara}) each term of (\ref{312bb}) is well-defined.
 Then, we can  estimate,  using (\ref{nvest}) and (\ref{vest})
 \bb\label{313bb}
|\nabla V(y,s)|\psi_{R_0}(y) \leq  |y|^\delta M \psi_{R_0} \leq R_0^\delta M \psi_{R_0} \leq  \frac{\a-\gamma}{4(\a+1)}\psi_{R_0} (y),
 \ee
\bq\label{314bb}
 |(V\cdot \nabla )\psi_{R_0}|&\leq& - |V|\partial_r \psi_{R_0} \leq  -M|y|^{1+\delta} \partial_r\psi_{R_0} \leq -MR_0^\delta|y|\partial_r\psi_{R_0}\n \\
 &\leq& -\frac{1}{2(\a+1)}  |y| \partial_r \psi_{R_0}
 =  -\frac{1}{2(\a+1)} (y\cdot \nabla)\psi_{R_0},
\eq
and
\bq\label{314abb}
\lefteqn{\frac{2|\gamma|}{\a+1} |y|^{2\gamma-1}  \left|V\cdot \frac{y}{|y|} \right|\psi_{R_0}\leq \frac{2|\gamma|}{\a+1}M|y|^{2\gamma-1}|y|^{1+\delta} \psi_R}\hspace{.3in} \n \\
&&\leq\frac{2|\gamma|}{\a+1}R_0^\delta M|y|^{2\gamma} \psi_{R_0} \leq
\frac{\a-\gamma}{2(\a+1)} |y|^{2\gamma} \psi_{R_0}
\eq
for all $s\in (0, S_0)$.
Substituting (\ref{313bb})-(\ref{314abb})  into (\ref{312bb}),  we  have
 \bqn
\lefteqn{\frac{\partial}{\partial s} (\psi_{R_0}|y|^{2\gamma}|B|^2 ) + \frac{(\a-\gamma)}{\a+1}\psi_{R_0}|y|^{2\gamma} |B|^2 +\frac{1}{\a+1}(y\cdot \nabla )( \psi_{R_0}|y|^{2\gamma}|B|^2 ) )}\n \\
 &&+(V\cdot \nabla )(\psi_{R_0}|y|^{2\gamma}|B|^2) \leq  \frac{1}{2(\a+1)} |y|^{2\gamma}|B|^2 (y\cdot \nabla )\psi_{R_0}\leq 0,
 \eqn
 and
 \bb\label{315a}
\frac{\partial f}{\partial s}  + \frac{(\a-\gamma)}{\a+1}f +\frac{1}{\a+1}(y\cdot \nabla )f
 +V\cdot \nabla f \leq 0,\quad f(y,s):=\psi_{R_0}|y|^{2\gamma}|B|^2.
 \ee
 Let $q>2$. Multiplying (\ref{315a}) by $f|f|^{q-2}$, and integrating over $B_{R_0}\times (0, S_0)$, and taking into account the periodicity, $f(y,0)=f(y, S_0)$, we obtain
$$ \left(\frac{\a-\gamma}{\a+1}-\frac{3}{q(\a+1)}\right)\int_0 ^{S_0} \int_{B_{R_0}} |f|^q dyds \leq 0.
$$
Choosing $q> \frac{3}{\a-\gamma}$, we obtain  $f=0$, and $B=0$ on the domain  $B_{R_0}\times (0, S_0)$. Therefore, the system (\ref{mhd}) reduces to the system $\se$ on $B_{R_0}\times (0, S_0)$. Repeating the proof of Theorem 1.2, we conclude
$V=0$ also on $B_{R_0}\times (0, S_0)$.\\

Let us choose  and $R\in (R_0, \infty)$ be fixed. We define here
\bb\label{rad2}
R_{k}=R_{k-1}+\mu^{\frac{1}{\delta}}, \quad k=1, \cdots, m
 \ee
where $M$ is the same as in (\ref{defm}), and $m$ is the smallest integer such that $\frac{R-R_0}{\mu^{\frac{1}{\delta}}}<m$.
 We  multiply  the second equations of (\ref{mhd}) by $B\psi_{R_k}|y|^{2\gamma} $ to obtain
\bq\label{312}
\lefteqn{\frac{\partial}{\partial s} (\psi_{R_k}|y|^{2\gamma}|B|^2 ) + \frac{2(\a-\gamma)}{\a+1}\psi_{R_k}|y|^{2\gamma} |B|^2 +\frac{1}{\a+1}(y\cdot \nabla )( \psi_{R_k}|y|^{2\gamma}|B|^2 ) )}\n \\
 &&+(V\cdot \nabla )(\psi_{R_k}|y|^{2\gamma}|B|^2) \leq 2|y|^{2\gamma}\psi_{R_k} |\nabla V| |B|^2+ \frac{1}{\a+1} |y|^{2\gamma}|B|^2 (y\cdot \nabla )\psi_{R_k} \n \\
 &&\hspace{1.in}\qquad\quad\qquad+|y|^{2\gamma}|B|^2 (V\cdot\nabla) \psi_{R_k} +\frac{2\gamma}{\a+1} |y|^{2\gamma-1} |B|^2 V\cdot\frac{y}{|y|}\psi_{R_k}.\n \\
 \eq
 We suppose it is shown that $V=B=0$ on $B_{R_{k-1}} \times (0, S_0)$, which is the case for $k=1$.   For $y\in B_{R_k}$ let us set $\bar{y}=\frac{y}{|y|} R_{k-1}$.  Then, we   estimate
 \bb\label{313}
|\nabla V(y,s)|\psi_{R_k}(y) \leq  |y-\bar{y}|^\delta M \psi_{R_k} \leq (R_k-R_{k-1})^\delta M \psi_{R_k} \leq  \frac{\a-\gamma}{4(\a+1)}\psi_{R_k} (y),
 \ee
\bq\label{314}
 |(V\cdot \nabla )\psi_{R_k}|&\leq& - |V|\partial_r \psi_{R_k} \leq  -M|y-\bar{y}|^{1+\delta} \partial_r\psi_{R_k} \leq -M(R_k-R_{k-1})^\delta|y|\partial_r\psi_{R_k}\n \\
 &\leq& -\frac{1}{2(\a+1)}  |y| \partial_r \psi_{R_k}
 =  -\frac{1}{2(\a+1)} (y\cdot \nabla)\psi_{R_k},
\eq
and
\bq\label{314a}
\lefteqn{\frac{2|\gamma|}{\a+1} |y|^{2\gamma-1}  \left|V\cdot \frac{y}{|y|} \right|\psi_{R_k}\leq \frac{2|\gamma|}{\a+1}M|y|^{2\gamma-1}|y-\bar{y}|^{1+\delta} \psi_R}\hspace{.3in} \n \\
&&\leq\frac{2|\gamma|}{\a+1}(R_k-R_{k-1})^\delta M|y|^{2\gamma} \psi_{R_k} \leq
\frac{\a-\gamma}{2(\a+1)} |y|^{2\gamma} \psi_{R_k}
\eq
for all $s\in (0, S_0)$.
Substituting (\ref{313})-(\ref{314a})  into (\ref{312}),  we  have
 \bqn
\lefteqn{\frac{\partial}{\partial s} (\psi_{R_k}|y|^{2\gamma}|B|^2 ) + \frac{(\a-\gamma)}{\a+1}\psi_{R_k}|y|^{2\gamma} |B|^2 +\frac{1}{\a+1}(y\cdot \nabla )( \psi_{R_k}|y|^{2\gamma}|B|^2 ) )}\n \\
 &&+(V\cdot \nabla )(\psi_{R_k}|y|^{2\gamma}|B|^2) \leq  \frac{1}{2(\a+1)} |y|^{2\gamma}|B|^2 (y\cdot \nabla )\psi_{R_k}\leq 0,
 \eqn
 and
 \bb\label{315a}
\frac{\partial f}{\partial s}  + \frac{(\a-\gamma)}{\a+1}f +\frac{1}{\a+1}(y\cdot \nabla )f
 +V\cdot \nabla f \leq 0,\quad f(y,s):=\psi_{R_k}|y|^{2\gamma}|B|^2.
 \ee
 Let $q>2$. Multiplying (\ref{315a}) by $f|f|^{q-2}$, and integrating over $B_{R_k}\times (0, S_0)$, and taking into account the periodicity, $f(y,0)=f(y, S_0)$, we obtain
$$ \left(\frac{\a-\gamma}{\a+1}-\frac{3}{q(\a+1)}\right)\int_0 ^{S_0} \int_{B_{R_k}} |f|^q dyds \leq 0.
$$
Choosing $q> \frac{3}{\a-\gamma}$, we obtain  $f=0$, and $B=0$ on the domain  $B_{R_k}\times (0, S_0)$. Therefore, the system (\ref{mhd}) reduces to the system $\se$ on $B_{R_k}\times (0, S_0)$, and $V=0$ on $B_{R_{k-1}}\times (0, S_0)$. By the proof of 1.2 we can extend the zero set of $V$  to $B_{R_k}\times (0, S_0)$.
Thus we have shown that $V=B=0$ on $B_{R_k}\times (0, S_0)$, and we can repeat the argument for $k=1, \cdots, m$, where $m$ is the smallest integer such that $\frac{R-R_0}{\mu^{\frac{1}{\delta}}}\leq m$, until we have the zero set of $V$ and $B$ as $B_{R-\mu^\frac{1}{\delta}}\times (0, S_0)$. Since $R$ can be arbitrary large, we conclude the proof.\\

  \noindent{\em \underline{(C) The case  with $\a <-1$:} }
  In this case we define $\overline{V} (y, s)= V(y, S_0-s)$,  $\bar{P} (y, s)= P(y, S_0-s)$ and $\overline{B} (y, s)= B(y, S_0-s)$ for $0\leq s\leq S_0$. Then, the second equations of (\ref{mhd}) become
 \bb\label{211}
\frac {\partial \overline{B}}{\partial s}-\frac{\a}{\a+1}\overline{B} - \frac{1}{\a+1}(y\cdot \nabla) \overline{B}-(\vb\cdot \nabla )\overline{B} =-(\overline{B}\cdot \nabla )\overline{V}.
\ee
 This is the same situation as the proof (B) above, since  $-\frac{\a}{\a+1}<0, -\frac{1}{\a+1} >0$. In particular, we note that the signs in front of the terms
 $(V\cdot \nabla )\O$ and $ (\O\cdot \nabla )V$ are not important in the estimates (\ref{313})-(\ref{314a}).
 Repeating the argument of (B) word by word we conclude $V=B=0$ on $\Bbb R^3\times (0, S_0)$. $\square$\\
\ \\
 $$ \mbox{\bf Acknowledgements } $$
 This research is supported partially by NRF Grants no.
 2006-0093854 and  no. 2009-0083521.

\end{document}